\documentclass[12pt]{article}
\usepackage{kz_modified}
\usepackage{enumerate}
\usepackage{authblk}
\usepackage[font=scriptsize,labelfont=bf]{caption}
\RequirePackage[colorlinks,linkcolor=blue,citecolor=blue,urlcolor=blue]{hyperref}
\newcommand{\blind}{1}

\begin{document}

\def\spacingset#1{\renewcommand{\baselinestretch}%
{#1}\small\normalsize} \spacingset{1}

%%%%%%%%%%%%%%%%%%%%%%%%%%%%%%%%%%%%%%%%%%%%%%%%%%%%%%%%%%%%%%%%%%%%%%%%%%%%%%

\if1\blind
{
  \title{\textbf{Discussion of `Multiscale Fisher's Independence Test for Multivariate Dependence'}}
  \author[1]{Duyeol Lee}
  \author[2]{Helal El-Zaatari}
  \author[2,3]{Michael R. Kosorok}
  \author[4]{Xinyi Li}
  \author[3]{Kai Zhang} 
  \affil[1]{Corporate Risk, Wells Fargo, McLean, VA 22102}    \affil[2]{Department of Biostatistics, University of North Carolina, Chapel Hill, NC 27599}
  \affil[3]{Department of Statistics and Operations Research, University of North Carolina, Chapel Hill, NC 27599}
\affil[4]{School of Mathematical and Statistical Sciences, Clemson University, Clemson, SC 29634}
  \maketitle
} \fi

\if0\blind
{
  \bigskip
  \bigskip
  \bigskip
  \begin{center}
    {\LARGE\bf BET on Independence}
\end{center}
  \medskip
} \fi

\noindent%

\spacingset{1.68} % DON'T change the spacing!

\section{Introduction}
The multiscale Fisher’s independence test (\textsc{MultiFIT} hereafter) proposed by \cite{10.1093/biomet/asac013} is a novel method to test independence between two random vectors. By its design, this test is particularly useful in detecting local dependence. Moreover, by adopting a resampling-free approach, it can easily accommodate massive sample sizes. Another benefit of the proposed method is its ability to interpret the nature of dependency. We congratulate the authors, Shai Gorksy and Li Ma, for their very interesting and elegant work. In this comment, we would like to discuss a general framework unifying the \textsc{MultiFIT} and other tests and compare it with the binary expansion randomized ensemble test (BERET hereafter) proposed by \cite{lee2019testing}. We also would like to contribute our thoughts on potential extensions of the method.

\section{\textsc{MultiFIT} under the binary expansion framework}

To understand the properties of \textsc{MultiFIT}, we analyze the test statistic under a similar multi-resolution approach in the binary expansion testing (BET hereafter) framework in \cite{zhang2019bet}. Consider the test of independence of two continuous variables from the copula $(U,V)$. The uniform consistency with respect to the total variation distance requires consistency for any alternative that is some distance from independence, i.e., $H_0: P_{(U,V)} = Unif[-1,1]^2 ~~\text{versus}~~ H_1: TV(P_{(U,V)},Unif[-1,1]^2) \ge \delta,~~\text{for some}~~ 0<\delta\le 1.$ Theorem~2.2 in \cite{zhang2019bet} shows the \textit{non-existence} of a test that is uniformly consistent with respect to the total variation distance. The key reason for non-uniform consistency is the unidentifiability in this problem. To avoid this issue and develop a nonparametric test statistic that is both powerful and robust, \cite{zhang2019bet} develops the BET framework to test approximate independence through a filtration approach. The filtration is constructed through the classical probability result of binary expansion: $U=\sum_{d=1}^\infty {A_{1,d}/2^d}$, $V=\sum_{d=1}^\infty {A_{2,d}/2^d}$ where $A_{1,d}$ and $A_{2,d}$ represent the $d$-th bit of $U$ and $V$ respectively. It is known that $U$ (or $V$) is marginally distributed as $Unif[-1,1],$ if and only if  $A_{1,d}\  ({\rm or}\ A_{2,d}) \stackrel{i.i.d.}{\sim} Rademacher$. When we truncate the expansions at finite depths $d_1$ and $d_2$,  $U_{d_1}=\sum_{d=1}^{d_1} {A_{1,d}/2^d}$ and $V_{d_2}=\sum_{d=1}^{d_2} {A_{2,d}/2^d}$ become discrete uniform variables that generate an analytically attractable filtration to approximate $(U,V)$. At every depth of the binary expansion filtration, the probability model of $(U_{d_1},V_{d_2})$ is a low resolution approximation of $(U,V)$ providing global distributional information. The interactions of these binary variables are denoted by $A_\Lambda$'s, where $\Lambda$ is a binary vector index with 1 or 0 indicating the presence of the bit of a variable in the interaction. The sums of observed interactions of $A_\Lambda$'s with $\Lambda \neq 0$ from $n$ samples, which are referred to as \textit{symmetry statistics} and are denoted by $\bar{S}_\Lambda=\sum_{i=1}^{n}A_{\Lambda,i}/n$, are shown to be completely sufficient for dependence and form the building blocks of inference.

\cite{zhang2021beauty} unifies several important tests such as the $\chi^2$ test, Spearman's $\rho$, and distance correlation in the binary expansion framework. It is shown that each of these test statistics can be approximated by a quadratic form of symmetry statistics $\bar{S}^TW \bar{S}$ for some deterministic weight matrix $W$ depending on the choice of distance. In the \textsc{MultiFIT} procedure, each of the four blocks in a cuboid corresponds to a cell in the contingency table from the discretization of the binary expansion approximation. The count of observations within each block can be written as a linear combination of symmetry statistics. For example, the test statistics for the first three cuboids corresponding to the case $p=2$, $k_1=1$ and $k_2=0$ are $(\sum_{i=1}^n A_{1,1,i}A_{2,1,i})^2, \sum_{i=1}^n (A_{1,2,i}A_{2,1,i}-A_{1,1,i}A_{1,2,i}A_{2,1,i})^2/2$ and $\sum_{i=1}^n (A_{1,2,i}A_{2,1,i}+A_{1,1,i}A_{1,2,i}A_{2,1,i})^2/2.$ In general, the test statistic $T_{R_{max}}$ of the \textsc{MultiFIT} procedure with the Bonferroni correction can be written in the following form as the maximum of many quadratic forms of symmetry statistics in the BET framework: $T_{R_{max}}=max_{1\leq j \le g(R_{max})}\{\bar{S}^T W_j\bar{S}\}$,
where $g(R_{max})$ is the number of cuboids as a function of $R_{max}$ and $W_j$ is a rank-one deterministic symmetric matrix corresponding to the $j$-th cuboid. With this binary expansion representation, the properties of \textsc{MultiFIT} can be investigated using the weights $W_1,\ldots,W_{g(R_{max})}$.

\section{Study of power}

The weight $W$ not only determines the form of the test but also determines its power properties. In general, if under the alternative the mean $\mu$ of the symmetry statistics vector falls in the eigenspace of the weight matrix corresponding to the highest (lowest) eigenvalues, then the test will have a high (low) power. Therefore, the determinateness of $W$ creates a key issue on the uniformity and robustness of the power, as it will always favor some alternatives but not others. In particular, the $W$'s of \textsc{MultiFIT} favor local dependencies rather than global ones. Approaches to avoid the above uniformity issue on the robustness of power include using some data-adaptive weights as in \cite{zhang2021beauty} and using an ensemble method as in \cite{lee2019testing}. In this comment, we focus on the comparison with BERET proposed by \cite{lee2019testing} which is also both robust and interpretable. We conduct simulation studies over the same simulation settings described in \cite{10.1093/biomet/asac013}. Figure~1 compares the performance of the methods. As discussed in \cite{10.1093/biomet/asac013}, \textsc{MultiFIT} has a natural advantage to detect marginal dependencies as it focuses on the testing of pairs of margins. As a result, it shows relatively better performance in marginal scenarios than in spread scenarios. While BERET shows robust performance across different scenarios, \textsc{MultiFIT} particularly outperforms in the local dependencies as expected. Because most existing tests are not designed to particularly test local dependencies, this power property can be a unique strength of \textsc{MultiFIT}. Additional investigations of the weight $W$ of \textsc{MultiFIT} are needed to more thoroughly reveal its usefulness in local dependencies.

\begin{figure}[htbp!]
\begin{center}
\begin{tabular}{cc}
	\includegraphics[scale=0.7]{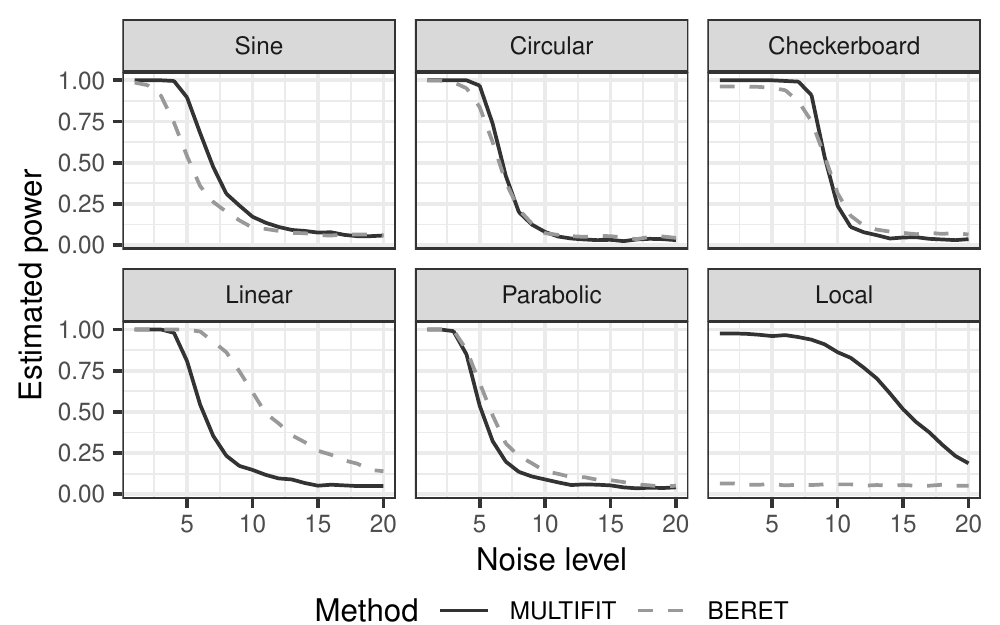}
	&\includegraphics[scale=0.7]{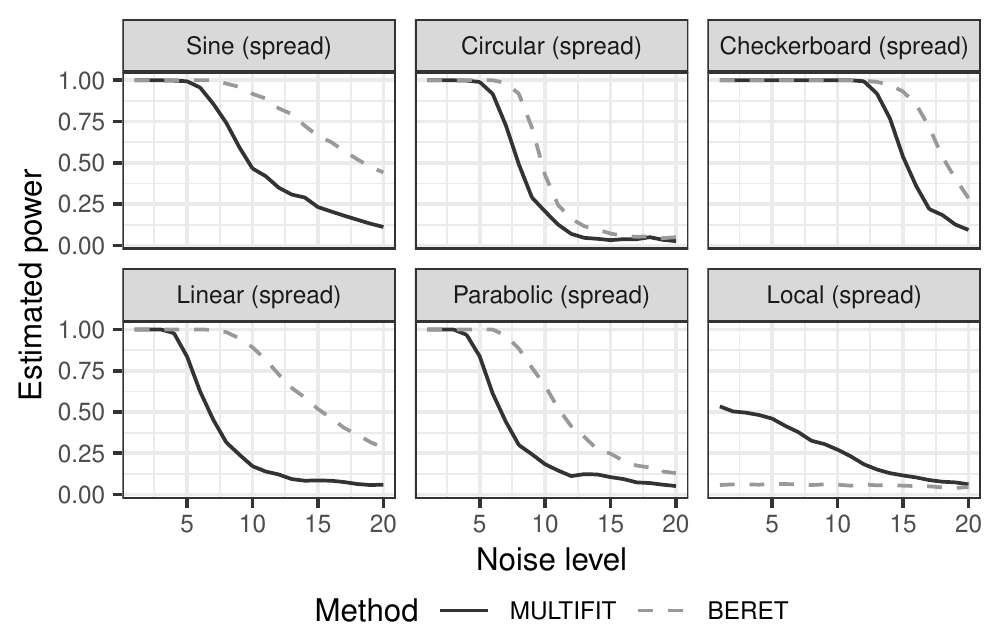}\\
	{(a)} & {(b)} \\
\end{tabular} 
\label{simresult}
\end{center}
\caption{Power versus noise level for different methods. (a) Estimated power at 20 noise levels for the different methods under the six marginal scenarios: \textsc{MultiFIT} (black solid), BERET with $d_{max}=4$ (grey dotted). (b) Estimated power at 20 noise levels under the six spread scenarios.}
\end{figure}

\section{Interpretability}

Clear interpretability is particularly important in evaluating multivariate dependencies and is another feature of \textsc{MultiFIT}. Both \textsc{MultiFIT} and BERET provide clear interpretability through the approach of multiple testing. The main difference between the two methods is their focuses on different dependencies. To compare the strengths of interpretability of the two methods, we use the rotated 3D circle example in \cite{10.1093/biomet/asac013}. This example generates a circular dependency between $X_3$ and $Y_3$ and then rotates it by $\pi/4$ degrees to hide the visible relationship. Figure~2 shows the strongest dependencies detected by the methods. \textsc{MultiFIT} returns the most significant local dependency within the local cuboid in the lower right corner, while BERET detects a global dependency contrasting the points between the white and dark regions with respect to a binary interaction.

\begin{figure}[htbp!]
\begin{center}
\begin{tabular}[t]{cccc}
	\includegraphics[scale=0.25]{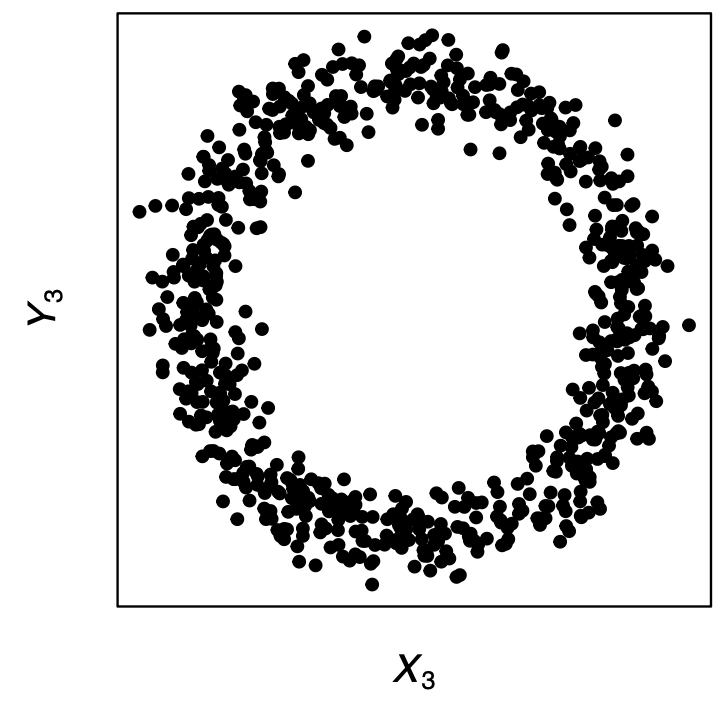}
	&\includegraphics[scale=0.25]{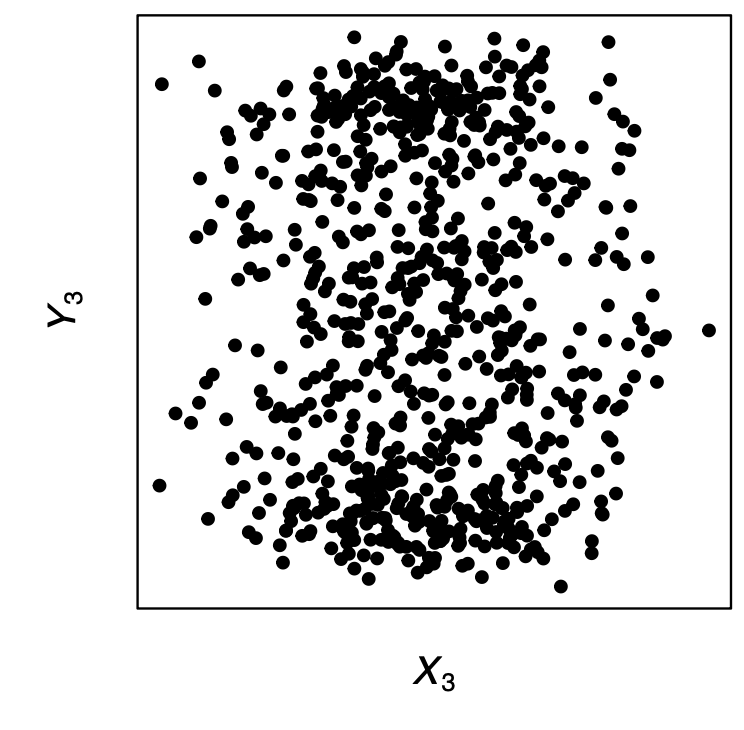}
	&\includegraphics[scale=0.25]{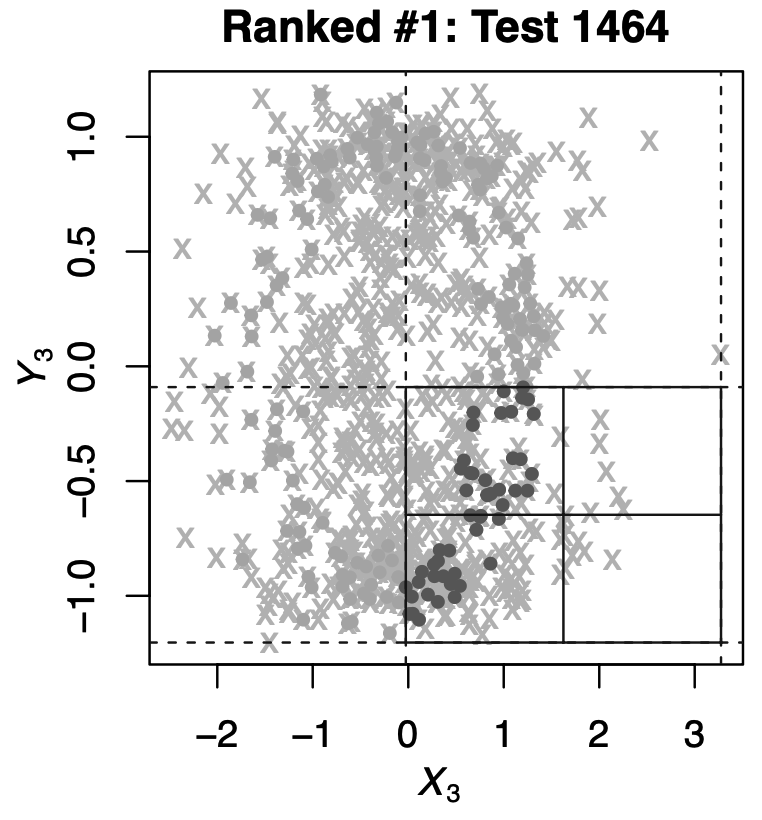}
	&\includegraphics[scale=0.46]{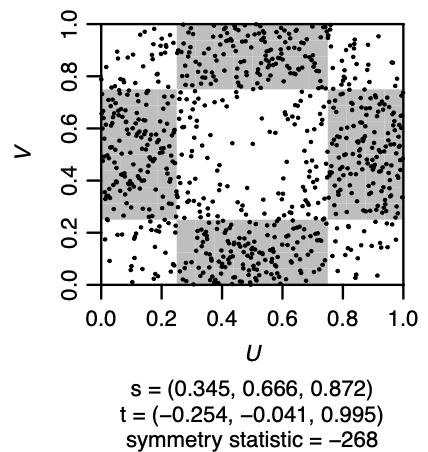}\\
	{(a)} & {(b)} & {(c)} & {(d)}\\
\end{tabular} 
\label{interpret}
\end{center}
\caption{(a) and (b) Marginal views of the 3D circle data before and after rotation. (c)  Scatter plot for the observations in the 2 × 2 tables identified by \textsc{MultiFIT}. (d) The strongest asymmetry detected by BERET with corresponding projections $s, t$ of $X$ and $Y$ and the symmetry statistic.}
\end{figure}

\section{Infinite dimensional extensions}
One potential extension is to infinite dimensional random variables, such as arises with functional data. For instance, the Brownian distance covariance statistic of \cite{szekely2009brownian} was shown in \cite{lyons2013distance} to be directly extendable to separable Hilbert spaces, with the same fundamental asymptotic properties as when applied to Euclidean random variables, provided the associated Hilbert space norm is used in the statistic instead of the Euclidean norm \citep{lyons2013distance} and the first moment of that norm is bounded. For example, suppose we want to assess the dependence between a random variable $X$, in a separable Hilbert space $H$ with norm $\|\cdot\|_H$, and a finite-dimensional Euclidean random vector $Y$. Then Brownian distance covariance would work provided $E(\|X\|_H)<\infty$ and $E(\|Y\|)<\infty$. Recall that a Hilbert space is separable if and only if it has a countable basis. An example of such a space is $L_2[0,1]$, and so functional data on a time interval rescaled to $[0,1]$ will work, along with many other richer Hilbert-space valued random variables. 

A question is how \textsc{MultiFIT} could similarly be extended to permit incorporation of Hilbert space valued random variables. One potential way forward would be applicable to Hilbert space random variables with finite second moment, $E(\|X\|_H^2)<\infty$, where we assume for simplicity of exposition that $X\in H=L_2[0,1]$ almost surely. In this situation, $X$ has a Karhunen-Lo{\'{e}}ve expansion of the form $X(t)=\nu(t)+\sum_{j=1}^{\infty}Z_j\lambda_j\phi_j(t)$, where $\nu(t)=E\{X(t)\}$, $Z_1,Z_2,\ldots$ are mean zero, variance 1, and mutually uncorrelated, and $(\lambda_j,\phi_j)$ are the principle component scores (i.e., eigenvalues) and associated eigenfunctions, and where $\phi_j$'s are an orthonormal basis for some subspace $H_0\subset H$ for which $\mathrm{pr}(X\in H_0)=1$ and the $\lambda_j$s are nonincreasing in $j$. Let $\tilde{X}_k=(Z_1,\ldots,Z_k)^T$ and let $X_k(t)$ be the Karhunen-Lo{\'{e}}ve expansion for $X$ summed up to the $k$th term. The map $M_k:\mathbb{R}^k\mapsto H$ which takes $\tilde{X}_k$ to $X_k$ is linear and, moreover, $X_k\rightarrow X$, in probability, as $k\rightarrow\infty$, since $\sum_{j=1}^{\infty}\lambda_j^2<\infty$ by the assumed existence of the second moment. This means that any stochastic dependence between $X$ and another random variable $Y$ (which is either Euclidean or separable Hilbert) will imply stochastic dependence between $X_k$ and $Y$ for some finite $k$. This structure can potentially be used to construct estimated finite Euclidean projections of Hilbert random variables which could be applied to \textsc{MultiFIT} for testing.

\bibliographystyle{biometrika}
\bibliography{discussion_multifit}

\end{document}